\documentclass[reqno,12pt]{amsart}

\usepackage{graphicx}
\usepackage{xcolor}
\usepackage{amssymb}
\usepackage{amsmath}
\usepackage{hyperref}

\theoremstyle{plain}
\newtheorem{theorem}{Theorem}

\newtheorem{rem}{Remark}
\newtheorem{question}{Question}

\theoremstyle{definition}

\theoremstyle{remark}

\newtheorem*{remark}{Remark}
\newtheorem*{remarks}{Remarks}

\def\N{{\mathbb N}}
\def\Z{{\mathbb Z}}

\def\R{{\mathbb R}}

\newcommand{\sysp}{\mathrm{sysp}}
\newcommand{\cro}{\mathrm{cr}}

\title{Sparse curve systems have intermediate growth type}

\date{\today}

\author{S.~Baader}
\email{sebastian.baader@unibe.ch}

\author{J.~J\"org}
\email{jasmin.joerg@unibe.ch}

\author{D.~Kosanovi\'c}
\email{danica.kosanovic@unibe.ch}

\address{Mathematisches Institut, Universit\"at Bern, Sidlerstrasse 5, 3012 Bern, Switzerland}

\begin{document}

\begin{abstract}
A system of simple closed curves on a surface of genus~$g$ is said to be sparse if their average pairwise intersection number does not exceed one. We show that the maximal size of a sparse curve systems grows roughly like a function of type $c^{\sqrt{g}}$, with $c$ between $2$ and $81938$.
\end{abstract}

\maketitle

\section{Introduction}

The closed orientable surface $\Sigma_g$ of genus $g \geq 2$ fits precisely $3g-3$ pairwise disjoint non-isotopic simple closed curves. If we relax the disjointness condition, we immediately run into difficult counting problems. For example, the maximal size of a $1$-system on $\Sigma_g$ --- that is, of a system of pairwise non-isotopic simple closed curves with pairwise intersection number either $0$ or $1$ --- is not known precisely. Thanks to the recent work of Greene~\cite{G}, the order of growth of this maximum is $g^2$, with a multiplicative error bounded by $\log(g)$. At the time of writing, the latter was removed by Aougab and Gaster, thus providing an almost sharp estimate on the maximal size of 1-systems \cite{AG}.

In this note, we consider a probabilistic version of $1$-systems: we say that a finite system $\Gamma$ of pairwise non-isotopic simple closed curves is \emph{sparse} if their average pairwise intersection number is at most $1$. Equivalently, the \emph{crossing number} $\cro(\Gamma)$ of $\Gamma$, defined as the sum of the intersection numbers of all pairs of curves of $\Gamma$, satisfies $\cro(\Gamma) \leq \binom{|\Gamma|}{2}$.

It is a priori not clear whether the size of a sparse curve system is a bounded function of the genus, that is, whether
\[
	\sysp(g):=\max\{|\Gamma| \mid \text{$\Gamma$ is a sparse curve system on $\Sigma_g$} \}
\]
is well-defined.
The main result in~\cite{BBS} implies that all but finitely many systems of systoles associated with congruence lattices in $\text{SL}(2,\Z)$ are sparse. These systems have polynomial growth in $g$, in analogy to $1$-systems. The following result might therefore come as a surprise. 

\begin{theorem}
\label{sparse1}
	For all $g \geq 16$ we have
	\[
	\frac{1}{16} [2\sqrt{g}] 2^{\sqrt{g}} \leq \sysp(g) \leq 2g e^{\sqrt{128g}+6}.
	\]
\end{theorem}
We use $[\sqrt{g}]$ for the floor of $\sqrt{g}$, as usual. A rough simplification of the bounds of Theorem~\ref{sparse1} for large $g$ yields
\[
	2^{\sqrt{g}} \leq \sysp(g) \leq 81938^{\sqrt{g}}.
\]
This suggests the following question.

\begin{question}
	Does the limit $\displaystyle{\lim_{g \to \infty} \frac{\log(\sysp(g))}{\sqrt{g}}}$ exist? If so, what is its value?
\end{question}

Theorem~\ref{sparse1} is a special case of the following stronger result. Let $f \colon \N \to \R_{>0}$ be any function. A finite system $\Gamma$ of pairwise non-isotopic simple closed curves on $\Sigma_g$ is \emph{$f(g)$-sparse} if their average pairwise intersection number is at most $f(g)$. Equivalently: $\cro(\Gamma) \leq f(g)\binom{|\Gamma|}{2}$. Generalizing the above, we define
\[
	\sysp_f(g):=\max\big\{|\Gamma| \mid \text{$\Gamma$ is an $f(g)$-sparse curve system on $\Sigma_g$} \big\}.
\]

\begin{theorem}
\label{sparse2}
	For all $\alpha \in (-1,1]$ and for all $g \geq 4^{\frac{2}{1+\alpha}}$ we have
	\[\frac{1}{16}[2g^{\frac{1-\alpha}{2}}]2^{g^{\frac{1+\alpha}{2}}}
     \leq \sysp_{g^\alpha}(g) \leq 2ge^{\sqrt{128} g^{\frac{1+\alpha}{2}}+6}.
	\]
\end{theorem}

The cases $\alpha=0$ and $\alpha=1$ deserve a special mention. The former is precisely Theorem~\ref{sparse1}; the second states that there exist exponential families of curves with pairwise average intersection number $g$. As we will see in the third section, there is a collapse in the limit case $\alpha=-1$: the size of a sparse $1/g$-system is bounded above by a linear function in $g$. In particular, for all $\alpha \leq -1$, the maximal size of a $g^\alpha$-system is essentially the same as the maximal size of a $0$-system, which is $3g-3$. On the other side, we do not know the maximal size of $g^\alpha$-systems for $\alpha>1$.

\begin{question}
	What is the growth type of the function $\sysp_{g^\alpha}(g)$ for $\alpha>1$?
\end{question}

The proof of Theorem~\ref{sparse2} has two parts. The first one is a construction of sparse curve systems whose size matches the lower bound. The second one is a derivation of the upper bound, which turns out to be a direct consequence of the recent crossing number estimate by Hubard and Parlier~\cite{HP}. Interestingly, the latter is valid for families of general closed curves, not only for embedded ones. In particular, it is not derived from Mirzakhani's curve counting statistics~\cite{M}, but from a result by Buser~\cite{B} (Lemma~6.6.4).

These two steps are carried out in the following two sections. 

\subsection*{Acknowledgements} 
We thank Luca Studer for providing the idea of using necklace-shaped surfaces in the construction of large curve systems. 

D.K.\ was supported by the SNSF Ambizione Grant PZ00P2-223790.

\section{Lower bound}

\begin{figure}[!htbp]
    \centering
    \includegraphics[width=0.2\linewidth]{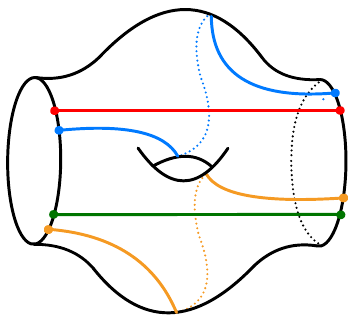}
    \caption{The four arcs on $\Sigma_{1,2}$.}
    \label{fig:4arcs}
\end{figure}

Consider the compact surface $\Sigma_{1,2}$  of genus 1 with 2 boundary components, and the four proper arcs on it as in Figure~\ref{fig:4arcs}. These arcs are simple, disjoint, pairwise non-isotopic, and non-boundary parallel. In fact, this is the maximal number of non-isotopic arcs connecting the two boundary components of $\Sigma_{1,2}$, by an Euler characteristic argument.

\begin{figure}[!htbp]
    \centering
    \includegraphics[width=0.86\linewidth]{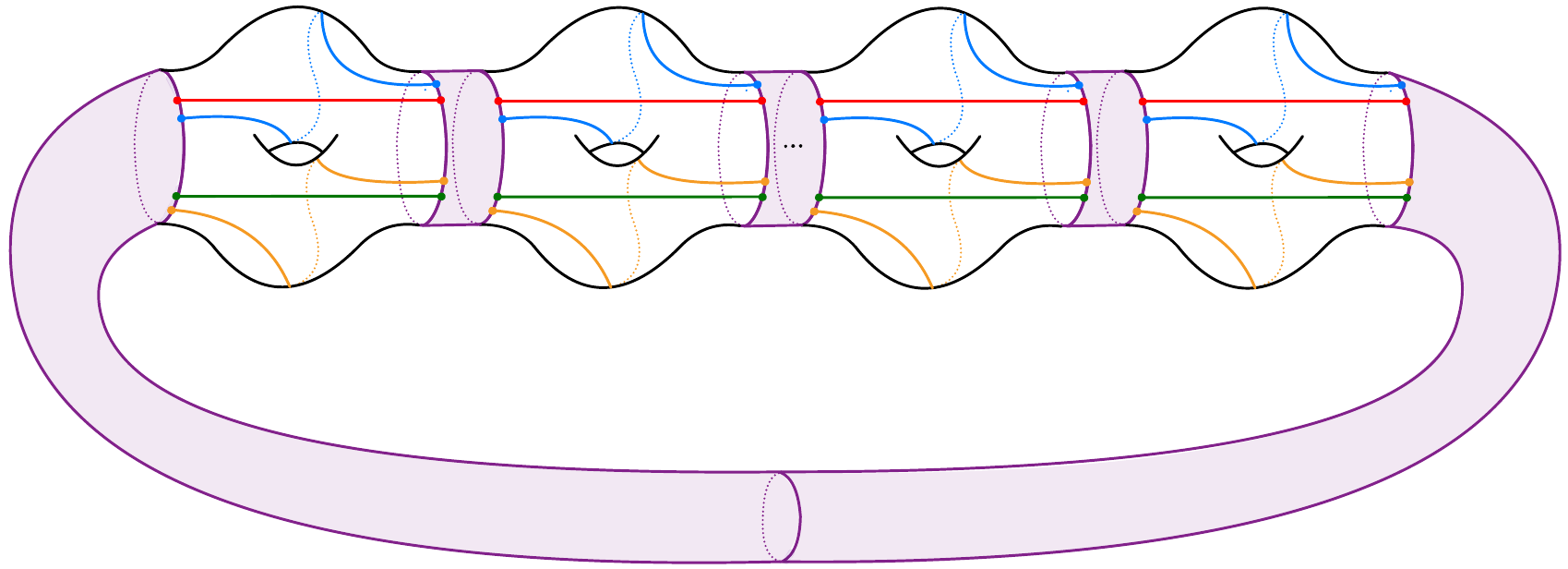}
    \caption{The genus $h$ surface $N$ consists of $h-1$ copies of $\Sigma_{1,2}$ and $h-1$ annuli.}
    \label{fig:necklace}
\end{figure}

Fix $\alpha\in(-1,1]$, and let $g\in\N$ satisfy $g\geq4^{\frac{2}{1+\alpha}}$. Then $h:=[\frac{1}{2}g^{\frac{1+\alpha}{2}}]$ satisfies $2\leq h\leq g$. Let us consider $h-1$ copies of $\Sigma_{1,2}$ and $h-1$ copies of the annulus $\Sigma_{0,2}$. We glue them together using a necklace shape as in Figure~\ref{fig:necklace}, thus forming a closed surface $N$ of genus $h$. 

On $N$ we now construct a collection $\Gamma_N$ of simple closed curves $\alpha_v$ for $v\in\{1,2,3,4\}^{h-1}$. In each copy of $\Sigma_{1,2}$ in $N$ we choose one of the four fixed arcs according to $v$, and in each annulus (shaded regions in Figure~\ref{fig:necklace}) we connect the ends of the respective arcs, by going to the right if necessary. No two of the curves $\alpha_v$ are homotopic since they represent different homology classes.

\begin{figure}[!htbp]
    \centering
    \includegraphics[width=0.35\linewidth]{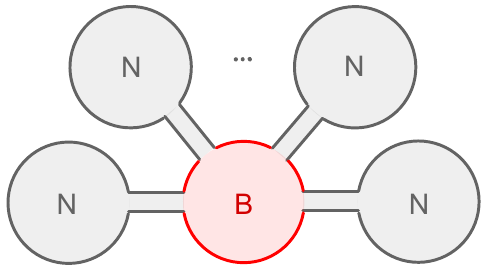}
    \caption{The surface $\Sigma_g$ consists of $h'$ copies of $N$ attached to the surface $B$.}
    \label{fig:schematic}
\end{figure}

Moreover, for any two of these $4^{h-1}$ curves in $\Gamma_N$, the crossings occur only in the $h-1$ shaded regions, and at most once therein. Hence,
\[
	\cro(\Gamma_N)\leq (h-1) \binom{4^{h-1}}{2}.
\]
Finally, we construct a curve system $\Gamma$ on a surface $\Sigma_g$ of genus $g$. We view $\Sigma_g$ as a base surface $B$ of genus $g-hh'$ to which we connect sum $h':=[2g^{\frac{1-\alpha}{2}}]$ copies of the necklace $N$; this is schematically depicted in Figure~\ref{fig:schematic}. 
Note that $hh'=[\frac{1}{2}g^{\frac{1+\alpha}{2}}][2g^{\frac{1-\alpha}{2}}]\leq \frac{1}{2}g^{\frac{1+\alpha}{2}}\cdot 2g^{\frac{1-\alpha}{2}}=g$, so $g-hh'$ is indeed a nonnegative integer.

For $\Gamma$, we take the union of the collections $\Gamma_N$ on each copy of $N$. Note that the only intersections that occur are those within each collection, as a curve never leaves its copy of $N$.

Thus, the total number of curves is 
\[
	|\Gamma|=h'\,4^{h-1}
\]
and their total number of crossings is 
\[
	\cro(\Gamma)\leq h'(h-1) \binom{4^{h-1}}{2}.
\]
Therefore, we have
\begin{align*}
	\frac{\cro(\Gamma)}{\binom{|\Gamma|}{2}}
	& \leq \frac
		{\frac{1}{2}h'(h-1) 4^{h-1} (4^{h-1}-1)}
		{\frac{1}{2}h'\,4^{h-1} (h'\,4^{h-1}-1)}
		\\
	& = \frac
		{(h-1)(4^{h-1} - 1)}
		{h'4^{h-1} - 1}
	=\frac
		{(h-1)(4^{h-1} - 1)}
		{h'(4^{h-1} - 1)+(h'-1)}
		\\
	& \leq \frac
		{(h-1)(4^{h-1}-1)}
		{h'(4^{h-1}-1)}
	= \frac{h-1}{h'}.
\end{align*}
Recall that we have defined $h=[\frac{1}{2}g^{\frac{1+\alpha}{2}}]$ and $h'=[2g^{\frac{1-\alpha}{2}}]$. Therefore, we can use the inequalities $h\leq \frac{1}{2}g^{\frac{1+\alpha}{2}}$ and $h'\geq 2g^{\frac{1-\alpha}{2}}-1$ to conclude 
\begin{align*}
	\frac{\cro(\Gamma)}{\binom{|\Gamma|}{2}} &\leq 
	\frac{h-1}{h'}
	\leq \frac
	{\frac{1}{2}g^{\frac{1+\alpha}{2}}-1}
	{2g^{\frac{1-\alpha}{2}}-1} \\
	&= \frac
    {g^{\frac{1+\alpha}{2}} (\frac{1}{2}-g^{-\frac{1+\alpha}{2}})}
    {g^{\frac{1-\alpha}{2}} (2-g^{-\frac{1-\alpha}{2}})}
    = g^\alpha \frac{\frac{1}{2}-g^{-\frac{1+\alpha}{2}}}{2-g^{-\frac{1-\alpha}{2}}} \\
    &\leq \frac{1}{2} g^\alpha.
\end{align*}
The last inequality follows by observing that the numerator is at most $\frac{1}{2}$ and the denominator is increasing in $g$ and thus bounded below by $2-1^{-\frac{1-\alpha}{2}}=1$.
Thus, $\Gamma$ is indeed a $g^\alpha$-sparse curve system with $[2g^{\frac{1-\alpha}{2}}]4^{[\frac{1}{2}g^{\frac{1+\alpha}{2}}]-1}$ curves on a genus $g$ surface. The lower bound in Theorem~\ref{sparse2} follows by observing that $[2g^{\frac{1-\alpha}{2}}]4^{[\frac{1}{2}g^{\frac{1+\alpha}{2}}]-1} \geq [2g^{\frac{1-\alpha}{2}}]4^{\frac{1}{2}g^{\frac{1+\alpha}{2}}-2} = \frac{1}{16}[2g^{\frac{1-\alpha}{2}}]2^{g^{\frac{1+\alpha}{2}}}$.

\begin{remark}
The above construction came out of the special case $\alpha=0$. In this case, both the number of subsurfaces $h'$ and their genus $h$ roughly coincide with $\sqrt{g}$.
\end{remark}

\section{Upper bound}

In this section, we deduce an upper bound for the size of an $f(g)$-sparse curve system from the crossing number inequality given by Hubard and Parlier \cite{HP}. In particular, for $f(g) = g^\alpha$, this yields the upper bound of Theorem~\ref{sparse2}.

Let $\Gamma$ be an $f(g)$-sparse curve system on a surface of genus $g\geq2$. Provided $|\Gamma| \geq e^6(2g-1)$, we may apply the crossing number inequality. By plugging $m=|\Gamma|$ into \cite[Theorem 1.2]{HP}, we obtain
\[
\frac{1}{128(2g-2)}\left(|\Gamma| \log\left(\frac{|\Gamma|}{(2g-1)e^6}\right)\right)^2 
< \cro(\Gamma).
\]
On the other hand, $\Gamma$ is $f(g)$-sparse, so $\cro(\Gamma)\leq f(g)\binom{|\Gamma|}{2}\leq  f(g) \frac{|\Gamma|^2}{2}$. Combining these two inequalities and rearranging the terms yields 
$\log\left(\frac{|\Gamma|}{(2g-1)e^6}\right)^2 < 64(2g-2)f(g)$, and therefore
\[ 
|\Gamma| < (2g-1)e^{\sqrt{64(2g-2)f(g)}+6}
< 2ge^{\sqrt{128g f(g)}+6}.
\]
For $|\Gamma|<e^6(2g-1)$ simply note that the bound is trivially satisfied: $e^6(2g-1) < 2ge^6 < 2ge^{\sqrt{128g f(g)}+6}$.

\begin{remarks}\quad
\begin{enumerate}
    \item For $f(g) = g^\alpha$, where $\alpha \leq -1$, this shows that a $g^\alpha$-sparse curve system grows at most linearly in $g$, thus exhibiting the same order of growth as a system of pairwise disjoint curves.
    \item The crossing number inequality applies more generally to closed (not necessarily simple) curves. The same upper bound, therefore, holds for $f(g)$-sparse systems of possibly non-simple curves. The fact that the same order of growth may be achieved by systems of simple curves stands in contrast to the asymptotic length distributions of simple and non-simple curves, respectively \cite{M,H}.
\end{enumerate}
\end{remarks}

\vspace{15pt}

\end{document}